\documentclass{gtmon_a}     % Basic GT/GTM/AGT style
\pdfoutput=1

\usepackage{amscd} 

\proceedingstitle{Proceedings of the School and Conference in Algebraic
Topology (The Vietnam National University, Hanoi, 9-20 August 2004)}
\conferencestart{9 August 2004}
\conferenceend{20 August 2004}
\conferencename{School and Conference in Algebraic Topology}
\conferencelocation{Vietnam National University, Hanoi, Vietnam}

\editor{John Hubbuck}
\givenname{John}
\surname{Hubbuck}

\editor{Nguy\~\ecircumflex{}n H V H\uhorn{}ng}
\givenname{H\uhorn{}ng}
\surname{Nguy\~\ecircumflex{}n}

\editor{Lionel Schwartz}
\givenname{Lionel}
\surname{Schwartz}

\title[The odd-primary Kudo--Araki--May algebra]{The odd-primary
Kudo--Araki--May algebra\vspace{-4pt}\\of algebraic Steenrod operations\\and invariant
theory}

\author{David J\,Pengelley}
\givenname{David J}
\surname{Pengelley}
\address{New Mexico State University\\\newline
Las Cruces NM 88003\\
USA}
\email{davidp@nmsu.edu}
\urladdr{}

\author{Frank Williams}
\givenname{Frank}
\surname{Williams}
\address{New Mexico State University\\\newline
Las Cruces NM 88003\\
USA}
\email{frank@nmsu.edu}
\urladdr{}

\volumenumber{11}
\issuenumber{}
\publicationyear{2007}
\papernumber{11}
\startpage{217}
\endpage{243}

\doi{}
\MR{}
\Zbl{}

\arxivreference{}


\received{1 March 2005}
\revised{24 January 2006}
\accepted{2 February 2006}
\published{14 November 2007}
\publishedonline{14 November 2007}
\proposed{}
\seconded{}
\corresponding{}
\editor{}
\version{}

\makeatletter
\def\cnewtheorem#1[#2]#3{\newtheorem{#1}{#3}[section]
\expandafter\let\csname c@#1\endcsname\c@theorem}

\def\S{Section }
\let\tiny\relax
\def\textquotedblleft{``}
\def\textquotedblright{''}

\theoremstyle{plain}
\newtheorem{theorem}{Theorem}[section]

\cnewtheorem{corollary}[theorem]{Corollary}

\cnewtheorem{conjecture}[theorem]{Conjecture}

\cnewtheorem{lemma}[theorem]{Lemma}
\cnewtheorem{proposition}[theorem]{Proposition}
\cnewtheorem{exercise}[theorem]{Exercise}
\theoremstyle{definition}
\cnewtheorem{notation}[theorem]{Notation}
\cnewtheorem{definition}[theorem]{Definition}
\cnewtheorem{example}[theorem]{Example}
\cnewtheorem{remark}[theorem]{Remark}

\makeatother  %  move after \newtheorem block

\numberwithin{equation}{section}

\begin{document}

\begin{asciiabstract}
We describe bialgebras of lower-indexed algebraic Steenrod operations over the
field with p elements, p an odd prime. These go beyond the operations
that can act nontrivially in topology, and their duals are closely related to
algebras of polynomial invariants under subgroups of the general linear groups
that contain the unipotent upper triangular groups. There are significant
differences between these algebras and the analogous one for p=2, in
particular in the nature and consequences of the defining Adem relations.
\end{asciiabstract}

\begin{htmlabstract}
We describe bialgebras of lower-indexed algebraic Steenrod operations over the
field with p elements, p an odd prime. These go beyond the operations
that can act nontrivially in topology, and their duals are closely related to
algebras of polynomial invariants under subgroups of the general linear groups
that contain the unipotent upper triangular groups. There are significant
differences between these algebras and the analogous one for p=2, in
particular in the nature and consequences of the defining Adem relations.
\end{htmlabstract}

\begin{abstract}   % type your abstract below
We describe bialgebras of lower-indexed algebraic Steenrod operations over the
field with $p$ elements, $p$ an odd prime. These go beyond the operations
that can act nontrivially in topology, and their duals are closely related to
algebras of polynomial invariants under subgroups of the general linear groups
that contain the unipotent upper triangular groups. There are significant
differences between these algebras and the analogous one for $p=2$, in
particular in the nature and consequences of the defining Adem relations.
\end{abstract}

\maketitle

\section{Introduction and statement of results}

Mod $p$ \textquotedblleft lower-indexed\textquotedblright\ operations $D_{i}$
$(i\geq0)$ arising via $\mathbb{F}_{p}$--equivariance, for odd primes $p$, were
constructed by Steenrod \cite{StEp} for the cohomology of topological spaces
and by Dyer and Lashof \cite{DyerLash} for the homology of iterated loop
spaces. The operations were described in a more general algebraic context by
May \cite{May}, who computed implicit relations among the $D_{i}$. Until
recently, however, all investigators immediately dropped the $D_{i}$ for $i$
not congruent to $0$ or $-1$ $\mathop{\rm mod} (  p-1 )  $, because
such operations act trivially in the cohomology and homology of topological
spaces \cite[page 104]{StEp}. Furthermore, for reasons of grading under
composition, those $D_{i}$ that can act nontrivially in topology were
converted to upper-indexed \textquotedblleft Steenrod\textquotedblright\ 
 operations, $P^{j}$ and $\beta P^{j}$ in cohomology, generating the
odd-primary Steenrod algebra with \textquotedblleft Adem\textquotedblright\ 
 relations, and $Q^{k}$ and $\beta Q^{k}$ in homology, generating the
Dyer--Lashof algebra. For instance, we recall that on even degrees in
cohomology, $\smash{D_{i}\co H^{2q} (  X;\mathbb{F}_{p} )  \rightarrow
H^{2qp-i} (  X;\mathbb{F}_{p} )} $, and for $\smash{u\in H^{2q} (
X;\mathbb{F}_{p} )}  $, Steenrod defined $\smash{P^{j}\co H^{2q} (
X;\mathbb{F}_{p} )  \rightarrow H^{2q+2 (  p-1 )  j} (
X;\mathbb{F}_{p} )}  $ by $\smash{P^{j}u= (  -1 )  ^{q-j}D_{2 (
p-1 )   (  q-j )  }u}$. This displays the discarding of most of
the operations and also shows that the composition algebras of operations
generated by the $D$'s versus the $P$'s will be dramatically different in
structure, especially since the degree of the underlying class, which varies
during composition, is involved in converting betwixt them.

The structure generated by all the $D_{i}$ subject to their universal
\textquotedblleft Adem\textquotedblright\ relations has a richness going
beyond topology and is finding application in the study of algebras of
polynomial invariants, of which we shall give a new example here. We have
also used this structure in work to appear \cite{PW3} to give a minimal
presentation for the mod $p$ cohomology of $\mathbb{C}P(\infty)$ as a module
over the Steenrod algebra, which in turn allows us to give a minimal
presentation of the cohomology of the classifying space $BU$ (ie the algebra
of symmetric invariants) as an algebra over the Steenrod algebra. Corresponding results at the prime $2$, some joint with Peterson, appear in
\cite{PPW,PW1,PW2}. In \cite{PW1} we analyzed the analogous algebra of
operations $D_{i}$ at the prime $2$ and named it the Kudo--Araki--May algebra
$\mathcal{K}$. Where results in this paper are completely analogous to those
in \cite{PW1}, we shall omit their proofs; proofs that are not analogous or
immediate will be given in subsequent sections.

A major contrast with odd primes is that at the prime $2$ all operations
$D_{i}$ can act nontrivially for spaces, and they all convert to Steenrod or
Dyer--Lashof operations, unlike the operations we are pursuing here. In work in
progress, of which we give an example in this paper, we extend the results of
Campbell \cite{Camp} and Kechagias \cite{Kech1994,Kech} to describe certain
algebras of odd primary polynomial invariants as subalgebras of the dual of
the bialgebra generated by the $D_{2i}$, along with their structure as
algebras over the Steenrod algebra. This will involve the broader set of
operations going beyond those that can act nontrivially for spaces.

We shall also see a surprising algebraic difference in the larger algebra of
operations generated by the $D_{2i}$: the Adem relations between these
operations are no longer entirely determined just by those generated by
\textquotedblleft inadmissibles\textquotedblright, as happened at the prime $2$ and in the odd-primary Steenrod algebra. Nonetheless our study of this
bialgebra and some of its subalgebras and quotient algebras, relying on
analysis of Adem relations via formal power series, will produce bases
consisting of certain (but not necessarily all) admissibles (see \fullref{defadmissible}), analogous to the Steenrod algebra, along with other
features such as a generalized Nishida action. Our Nishida action provides
structure over the Steenrod algebra that can be compared with the action
independently enjoyed by algebras of invariants.

In our applications we deal with polynomial invariants or cohomology of spaces
concentrated in even degrees, so we shall simplify by eliminating Bocksteins
and thus deal only with the even $D$'s, which we shall denote by $e_{i}%
=D_{2i}$, for $i\geq0.$ We begin our analysis of the algebraic structure
generated by the $e_{i}$, and adopt the notational convention henceforth that
$e_{i}=0$ unless $i$ is a nonnegative integer.

To prepare for applications to polynomial invariants, and because we must
study the implementation of the Adem relations very carefully, we begin
formally with just the free algebra on the $e_{i}$, before the imposition of
any Adem relations.

\begin{definition}
Let $\widehat{\mathcal{U}}$ be the free noncommutative $\mathbb{F}_{p}%
$--algebra generated by elements $e_{i}$, for $i\geq0$, of \textit{topological
\/} degree $\left\vert e_{i}\right\vert =2i.$ The topological degree of a
product $xy$ in $\smash{\widehat{\mathcal{U}}}$ is given by $\left\vert xy\right\vert
=\left\vert x\right\vert +p\left\vert y\right\vert $ (cf \cite{PW1}). Note
that $\smash{\widehat{\mathcal{U}}}$ is bigraded, by topological degree $t$ and by
length $n$. We write $\smash{\widehat{\mathcal{U}}_{n,t}}$ for the component in this
bidegree. We caution that the reader should not assume that $e_{0}$ is the
identity. It is not: $e_{0}\in\smash{\widehat{\mathcal{U}}_{1,0}}$, whereas
$1\in\smash{\widehat{\mathcal{U}}_{0,0}}$. Note too that the nature of the formula for
topological degree on products, in fact the principal purpose of its skewed
nature, will be to make all relations homogeneous in topological degree in
addition to length; this will be apparent from any of the formulations of the
relations below.
\end{definition}

If we restrict attention to May's algebraic Steenrod operations $e_{i}$
applied only to even-dimensional classes, we see that his universal formulas
\cite[page 180]{May} involving the algebraic Steenrod operations are induced by
imposing on $\widehat{\mathcal{U}}$ the relations
\begin{multline*}
  \sum_{k}(-1)^{k+s}\binom{s-(p-1)k}{k}e_{r+(pk-s)(p-1)}e_{s-k(p-1)}\\
  \sim\sum_{l}(-1)^{l+r}\binom{r-(p-1)l}{l}e_{s+(pl-r)(p-1)}e_{r-l(p-1)}%
\text{ }\\
 \text{for each fixed }r,s\geq0\text{,}%
\end{multline*}
which we shall refer to as May's relations.

We note that since $e_{i}=0$ for $i$ negative, these sums are finite for each
pair $r,s$. But these relations are extremely difficult to use in practice;
for instance, a particular monomial may appear in multiple relations on either
side. Happily, the effect of these elaborate relations can be unraveled into
equivalent relations that are more tractable and of more familiar form, first
by interpreting them in terms of formal power series by defining
\[
e(u)=\sum_{i=0}^{\infty}e_{i}u^{i}.
\]
A short calculation followed by a change of variables gives the following theorem.

\begin{theorem}
\label{efps}May's relations above for the $e_{i}$ can be encoded as a formal
power series identity expressing a certain symmetry:
\[
e(u)e(v(v^{p-1}-u^{p-1}))\sim e(v)e(u(u^{p-1}-v^{p-1})),
\]
in which the coefficients of the monomials $u^{r}v^{s}$ are the corresponding
individual relations above.
\end{theorem}

To obtain equivalent relations that have many of the useful features of the
familiar Adem relations for the Steenrod and Dyer--Lashof algebras, we can use
the residue method of Bullett and Macdonald \cite{BullMac} and Steiner
\cite{Stei}, as we did in \cite[page 1461]{PW1},
to obtain the next theorem.

\begin{theorem}
[proved in \fullref{sec-a}]\label{eademcoord}May's relations are
equivalent to the relations
\[
e_{i}e_{j}\sim\sum_{k}(-1)^{\frac{pk-i}{p-1}}\binom{k-j-1}{\frac{pk-i}{p-1}%
-j}e_{i+pj-pk}e_{k}\qquad\text{for all integers }i,j\in\mathbb{Z},
\]
where the numerator in a binomial coefficient may be any integer, and a term
is present only if the fraction shown is an integer.
\end{theorem}

The very complicated relations of May can thus be replaced by these equivalent
relations, which are at least somewhat like traditional Adem relations, in the
sense that each two-fold monomial is now related only to a single sum of
two-fold monomials. However, notice that these replacement relations can be
nontrivial for arbitrary integer indices $i$ and $j$, since even though the
left side is automatically zero if $i$ or $j$ is negative, the right side may
not be, for instance, $e_{-\left(  p-1\right)  }e_{p}\sim-2e_{1}e_{p-1}$.
Together we shall call these the \emph{full relations\/}, since they are
bi-indexed by all $i,j\in\mathbb{Z}$, whereas we will call the proper subset
that is bi-indexed only by $i,j\geq0$ the \emph{Adem relations\/}.

Since the monomials $e_{i}e_{j}$ for $i,j\geq0$ form a basis for the two-fold
operations before relations are imposed, we wish to do our linear algebra with
only the Adem relations, rather than the full relations. Thus we will
initially develop our resulting quotient algebra(s) formed by imposing only
the Adem relations. Once we understand various important aspects of these, we
will be in a position to prove, in \fullref{negredundant}, that the
additional relations, ie for negative $i$ or $j$, are actually redundant.

We shall also be interested in two sub/quotient algebras of $\widehat
{\mathcal{U}}$:

\begin{definition}
Let $\widetilde{\mathcal{U}}$ be the subalgebra of $\widehat{\mathcal{U}}$
generated by the $e_{i}$ for which $i$ is even, and let $\mathcal{U}$ be the
subalgebra of $\widehat{\mathcal{U}}$ generated by the $e_{i}$ for which $i$
is divisible by $p-1.$ In dealing with $\mathcal{U}$ we shall use the
notation $d_{i}=e_{i(p-1)}$ for its algebra generators, all brought together
in the formal power series notation $d(t)=\sum_{i=0}^{\infty}d_{i}t^{i}$. (On
spaces concentrated in even degrees, only the operations $d_{i}$ can act
nontrivially \cite[page 104]{StEp}. It was the $d_{i}$ that were originally used
by Steenrod to construct the reduced power operations $P^{j}$ of the Steenrod
algebra $\mathcal{A}$.) Note that $\widetilde{\mathcal{U}}$ and $\mathcal{U}$
can also be regarded as quotient algebras, namely as $\widehat{\mathcal{U}}$
modulo the two-sided ideals generated by the $e_{i}$ for which $i$ is odd, and
the $e_{i}$ for which $i$ is not divisible by $p-1$, respectively. In fact
this is how they will be considered henceforth.
\end{definition}

We give $\widehat{\mathcal{U}}$ the structure of a bialgebra by defining a
component coalgebra structure with diagonal map given by
\[
\Delta(e_{i})=\sum_{a}e_{a}\otimes e_{i-a}.
\]
In formal power series, this becomes
\[
\Delta e(u)=e(u)\otimes e(u),
\]
ie $e(u)$ is \emph{grouplike\/}. With this definition, we see that
$\widetilde{\mathcal{U}}$ and $\mathcal{U}$ become quotient bialgebras.

To prepare for our linear algebra on $\widehat{\mathcal{U}}$, $\widetilde
{\mathcal{U}}$, and $\mathcal{U}$, note that since there is exactly one Adem
relation for each basis element of $\widehat{\mathcal{U}}_{2,\ast}$, we can
encode them in an endomorphism.

\begin{definition}
Let $\theta\co\widehat{\mathcal{U}}_{2,\ast}$ $\rightarrow$ $\widehat
{\mathcal{U}}_{2,\ast}$ be defined by the formula
\[
\theta(e_{i}e_{j})=\sum_{k}(-1)^{\frac{pk-i}{p-1}}\binom{k-j-1}{\frac
{pk-i}{p-1}-j}e_{i+pj-pk}e_{k}\qquad\text{for all }i,j\geq0,
\]
in other words, assign the right side of an Adem relation to its left side for
each basis element $e_{i}e_{j}\in\widehat{\mathcal{U}}_{2,\ast}$. NB: Here it
is critical that $i,j\geq0$, since we are defining $\theta$ using a basis of
$\widehat{\mathcal{U}}_{2,\ast}$.
\end{definition}

Since we will sometimes also need to use the right side of a full relation
even if $i,j$ do not satisfy the nonnegativity requirement of the definition
of $\theta$ on basis elements of $\widehat{\mathcal{U}}_{2,\ast}$, we extend
the notation to all $i,j\in\mathbb{Z}$.

\begin{definition}
Let the notation $\theta(i,j)$ be defined by
\[
\theta(i,j)=\sum_{k}(-1)^{\frac{pk-i}{p-1}}\binom{k-j-1}{\frac{pk-i}{p-1}%
-j}e_{i+pj-pk}e_{k}\qquad\text{for all }i,j\in\mathbb{Z},
\]
recognizing that this is merely a function on pairs $i,j\in\mathbb{Z}$, since
$e_{i}e_{j}$ is a basis element of $\smash{\widehat{\mathcal{U}}_{2,\ast}}$ only when
$i,j\geq0.$ In fact $e_{i}e_{j}=0$ by definition when $i$ or $j$ is negative,
while $\theta(i,j)$ as defined may not be.
\end{definition}

Regarding relations induced on $\widetilde{\mathcal{U}}$ and $\mathcal{U}$,
careful examination of the full relations on $\widehat{\mathcal{U}}$ in
\fullref{eademcoord} shows that, on the quotients (or subalgebras)
$\widetilde{\mathcal{U}}$ and $\mathcal{U}$, they induce full relations, and
also a corresponding endomorphism $\theta$ using the Adem relations. We also
have notation $\theta(i,j)$ for all appropriate $i,j\in\mathbb{Z}$, merely by
restricting all subscripts to those of the respective generators of each
quotient. In particular the full relations on $\mathcal{U}$ can be written
more succinctly.

\begin{theorem}
\label{dademcoord}In $\mathcal{U}$ the full relations are
\[
d_{i}d_{j}\sim\sum_{l}\left(  -1\right)  ^{pl-i}\binom{\left(  p-1\right)
\left(  l-j\right)  -1}{pl-i-\left(  p-1\right)  j}d_{i+pj-pl}d_{l}%
\qquad\text{for all }i,j\in\mathbb{Z}.
\]
\end{theorem}

For expressing a formal power series identity in $\mathcal{U}$ equivalent to
these, and envisaging the analogy to what we began with for May's relations in
$\widehat{\mathcal{U}}$, we first define two helpful functions.

\begin{definition}
(1)\qua $\varphi(a,b)=a(a^{p-1}-b^{p-1})$, and

(2)\qua $\psi(a,b)=a(a-b)^{p-1}.$
\end{definition}

Note that the formal power series identity for May's relations in
$\widehat{\mathcal{U}}$ can then be written just as
\[
e(u)e(\varphi(v,u))\sim e(v)e(\varphi(u,v)).
\]
The full relations on the quotient $\mathcal{U}$ are then equivalent to the
following formal power series identity.

\begin{theorem}
[proved in \fullref{sec-a}]\label{dfps}The relations in $\mathcal{U}$
induced by May's relations in $\widehat{\mathcal{U}}$ are%
\begin{gather*}
d(u)d(\psi(v,u))\sim d(v)d(\psi(u,v))\text{, ie }\\
d(u)d(v(v-u)^{p-1})\sim d(v)d(u(u-v)^{p-1}).
\end{gather*}
\end{theorem}

Now we are ready to consider the quotients by Adem relations.

\begin{definition}
We shall denote the algebra quotients of $\widehat{\mathcal{U}}$,
$\widetilde{\mathcal{U}}$, and $\mathcal{U}$ by their Adem relations by
$\widehat{\mathcal{K}}$, $\widetilde{\mathcal{K}}$ and $\mathcal{K}$,
respectively. (Recall that later we will show that the full relations with
negative $i$ or $j$ are redundant.) Note that since, as remarked above, the
relations are homogeneous in both bidegrees, $\widehat{\mathcal{K}}$,
$\widetilde{\mathcal{K}}$, and $\mathcal{K}$ will inherit these bidegrees as
well. By analogy with the prime $2$ \cite{PW1}, and because it consists of the
algebra of operations acting in topology on even degrees, we shall call
$\mathcal{K}$ the even topological Kudo--Araki--May algebra. Obviously
$\widehat{\mathcal{K}}$ and $\widetilde{\mathcal{K}}$ provide larger, purely
algebraic, versions.
\end{definition}

Relationships between $\mathcal{K}$ and the Steenrod and Dyer--Lashof algebras
are inherent in the consequences of the conversion formulas \cite[pp 161--2,
182]{May}%
\begin{align*}
P^{j}u&=\left(  -1\right)  ^{q-j}d_{q-j}u,\text{ where }u\text{ is a cohomology
class of degree }2q,\\
\tag*{\hbox{and}}
Q^{j}u&=\left(  -1\right)  ^{j-q}d_{j-q}u,\text{ where }u\text{ is a homology
class of degree }2q.
\end{align*}
Under composition using these conversions, the Adem relations in $\mathcal{K}$
produce, respectively, the traditional Adem relations in the Steenrod and
Dyer--Lashof algebras. However, the reader should not imagine that this
conversion provides anything as simple as an algebra map between $\mathcal{K}$
and either of the other two algebras. For instance, in either case the
operation $d_{0}=e_{0}$ always converts to the unstable $p$--th power operation
on any class. Since the degree of a class is involved in the conversion of any
operation on that class, and this degree changes during composition, the
relationship that arises is that of a \textquotedblleft sheared algebra
map\textquotedblright\ as described in \cite{PW1}.

We are now ready to describe bases for $\widehat{\mathcal{K}}$, $\widetilde
{\mathcal{K}}$, and $\mathcal{K}$.

\begin{definition}
\label{defadmissible}We shall call a monomial $e_{i_{1}}\cdots e_{i_{n}}$ or
$d_{i_{1}}\cdots d_{i_{n}}$ in any of our algebras \textit{admissible\/} if
$i_{1}\leq\cdots\leq i_{n}$, otherwise \textit{inadmissible\/}. And if
$e_{i}e_{j}$ is admissible, we call the nonnegative number $j-i$ its
\textit{excess\/}.
\end{definition}

For $p=2$ the inadmissible Adem relations (ie for $i>j$) completely
determine the admissible relations (ie for $i\leq j$), making the latter
redundant \cite{Le,PW1}, and thus the admissible monomials form a basis. We
shall see that this is also the case at odd primes for the Adem relations in
$\mathcal{U}$ (with generators $d_{i}=e_{i\left(  p-1\right)  }$) that create
$\mathcal{K}$, and thus the admissible monomials form a basis for
$\mathcal{K}$ (and correspondingly for the Steenrod algebra) at odd
primes. However, the situation is very different in $\widehat{\mathcal{U}}$
and in $\widetilde{\mathcal{U}}$, where the subscripts on generators are not
necessarily divisible by $p-1$. The following proposition illustrates some of
the serious consequences of the nonredundancy of admissible Adem relations
outside $\mathcal{U}$ and $\mathcal{K}$.

\begin{proposition}
[proved in \fullref{sec-b}]\label{noncongruenttrivialproducts}

{\rm (a)}\qua In the quotient $\widehat{\mathcal{K}}$, if $i$ $\!\not\equiv\!$
$j$ $\operatorname{mod}\left(  p-1\right)  $ then $e_{i}e_{j}=0$.

{\rm (b)}\qua In the quotient $\widehat{\mathcal{K}}$, if either $i$ or $j$ is
odd, then $e_{i}e_{j}=0$.
\end{proposition}

Since part (b) tells us that any products not occurring entirely inside
$\widetilde{\mathcal{K}}$ are zero, we see that $\widehat{\mathcal{K}}$ and
$\widetilde{\mathcal{K}}$ are the same in lengths greater than one. And even
within $\widetilde{\mathcal{K}}$, part (a) tells us that admissibles not
having mutually congruent indices will be zero, forcing a sparseness to any
possible basis. We are able to show that these two phenomena are the total
extent of the collapsing effects of the admissible Adem relations, leaving the
rest intact.

\begin{theorem}
[proved in \fullref{sec-b}]\label{basis}

{\rm (a)}\qua A vector space basis for $\widetilde{\mathcal{K}}$ (and for
$\widehat{\mathcal{K}}$ in lengths exceeding one) is given by the monomials
$e_{i_{1}}\cdots e_{i_{n}}$ for which all $i_{k\text{ }}$ are even, $0\leq
i_{1}\leq\cdots\leq i_{n}$, and $i_{k}\equiv i_{k-1}$ $\operatorname{mod}%
\left(  p-1\right)  $ for $2\leq k\leq n.$ 

{\rm (b)}\qua A vector space basis for $\mathcal{K}$ is given by the monomials
$d_{i_{1}}\cdots d_{i_{n}}$ for which $0\leq i_{1}\leq\cdots\leq i_{n\text{ }%
}.$
\end{theorem}

At this point we introduce important endomorphisms, which are essential for
proving that the basis theorem above ensures that the negative relations are
redundant, as alluded to earlier, and which then produce endomorphisms of the
quotients by the Adem relations.

Define the three algebra maps $\widehat{\alpha}\co\widehat{\mathcal{U}%
}\rightarrow\widehat{\mathcal{U}}$, $\wtilde{\alpha}\co\widetilde
{\mathcal{U}}\rightarrow\widetilde{\mathcal{U}}$, and $\alpha\co\mathcal{U}%
\rightarrow\mathcal{U}$ by the formulas $\widehat{\alpha}(  e_{i})
=e_{i-1}$, $\wtilde{\alpha}(  e_{i})  =e_{i-2}$, $\alpha
(  e_{i})  =e_{i-(  p-1)  }$ on their respective algebra generators.

\begin{theorem}
[proved in \fullref{sec-b}]\label{negredundant}The basis elements of
$\widehat{\mathcal{K}}$, $\widetilde{\mathcal{K}}$ and $\mathcal{K}$ given by
\fullref{basis} remain linearly independent if we impose the additional
(negative) relations
\[
e_{i}e_{j}\sim\sum_{k}(-1)^{\frac{pk-i}{p-1}}\binom{k-j-1}{\frac{pk-i}{p-1}%
-j}e_{i+pj-pk}e_{k}\qquad\text{for }i<0\text{ or }j<0.
\]
Thus imposing either the full relations or just the Adem relations produces
the same quotients.
\end{theorem}

Now that we know that the full relations and their proper subset the Adem
relations are equivalent impositions, we may use them interchangeably in
studying the effects of imposing them. In particular, it is easy to see from
comparing their formal power series formulations in \fullref{efps} and
\fullref{dfps} with the power series formulations $\widehat{\alpha}\left(
e(u)\right)  =ue(u)$, $\wtilde{\alpha}\left(  e(u)\right)  =u^{2}e(u)$, and
$\alpha\left(  d(u)\right)  =ud(u)$, that $\wtilde{\alpha}$ and $\alpha$
commute with their respective relations (which can also be verified by direct
calculation). Thus they induce algebra endomorphisms on $\widetilde
{\mathcal{K}}$ and $\mathcal{K}$, respectively. Note, however, that
$\widehat{\alpha}$ fails to commute with its associated relations by a sign.

NB: Since $\widehat{\mathcal{K}}$ degenerates to $\widetilde{\mathcal{K}}$ in
lengths greater than two, we focus most of our attention henceforth on
$\widetilde{\mathcal{K}}$ and $\mathcal{K}$.

The next theorem follows immediately from the formal power series formulas for
the relations and for $\Delta.$

\begin{theorem}
The diagonal maps $\Delta$ in $\widetilde{\mathcal{U}}$ and $\mathcal{U}$
respect the relations, and hence $\widetilde{\mathcal{K}}$ and $\mathcal{K}$
inherit the structure of bialgebras.
\end{theorem}

In preparation for connections to algebras of invariants via dualization, we
are interested in the coalgebra primitives in these bialgebras. The components
of each of these bialgebras in fixed length degree $n$ are coalgebras. Their
primitive elements are given in the following theorem. Parts 1, 2, and 3 are
immediate, while parts 4 and 5 are proved analogously to Theorem A of
\cite{PW1}.

\begin{theorem}
{\rm (1)}\qua A basis for the coalgebra primitives in $\widehat{\mathcal{U}%
}_{n,\ast}$ consists of the elements $e_{0}^{n}$ and $e_{0}^{a}e_{1}%
e_{0}^{n-a-1}$, for $0\leq a\leq n-1.$ 

{\rm (2)}\qua For those in $\widetilde{\mathcal{U}}_{n,\ast}$: $e_{0}^{n}$ and
$e_{0}^{a}e_{2}e_{0}^{n-a-1}$, for $0\leq a\leq n-1.$ 

{\rm (3)}\qua For those in $\mathcal{U}_{n,\ast}$: $d_{0}^{n}=e_{0}^{n}$ and
$d_{0}^{a}d_{1}d_{0}^{n-a-1}=e_{0}^{a}e_{p-1}e_{0}^{n-a-1}$, for $0\leq a\leq
n-1.$ 

{\rm (4)}\qua For those in $\smash{\widetilde{\mathcal{K}}_{n,\ast}}$: $e_{2}^{n}$ and
$e_{0}^{a}e_{p-1}^{n-a}$, for $1\leq a\leq n.$ 

{\rm (5)}\qua For those in $\mathcal{K}_{n,\ast}$: $d_{0}^{a}d_{1}^{n-a}%
=e_{0}^{a}e_{p-1}^{n-a}$, for $0\leq a\leq n.$
\end{theorem}

Now we discuss applications via dualization connecting these coalgebras to
algebras of polynomial invariants over the Steenrod algebra $\mathcal{A}$.
There are already some known results of interest. For instance, it was proved
by Kechagias in \cite[Theorem 2.23]{Kech1994} that in length degree $n$, the dual
algebra to the coalgebra $\widehat{\mathcal{U}}_{n,\ast}$ is isomorphic to the
invariants of the polynomial ring $S=\smash{\mathbb{F}_{p}[t_{1},\ldots,t_{n}]}$ in
$n$ variables of degree $2$, under the action of $\smash{\widehat{T}_{n}}$, the upper
triangular group with $1$'s down the principal diagonal. We shall discuss more
examples below, but first we wish to introduce additional structure on our
coalgebras for comparison with the natural action of the Steenrod algebra on
algebras of invariants.

We shall create a downward action of $\mathcal{K}^{\text{op}}$ (and hence of
$\mathcal{A}^{\text{op}}$) on $\widehat{\mathcal{U}}$,
\[
\mathcal{K}_{m,i}^{\text{op}}\otimes\widehat{\mathcal{U}}_{n,j}\rightarrow
\widehat{\mathcal{U}}_{n,\upnfrac{i+j}{p^{m}}},
\]
whose contragradient then automatically produces an unstable action of
$\mathcal{A}$ on the dual of $\widehat{\mathcal{U}}$. In \fullref{slinv}
below, and in work in progress, we will relate this to the action of
$\mathcal{A}$ on certain algebras of invariants.

Analogously to the prime $2$ \cite{PW1}, we denote this action by the symbol
$\ast$, refer to it as the \textit{Nishida\/} action, since it generalizes the
interaction discovered by Nishida between the actions of the Steenrod and
Dyer--Lashof algebras on the homology of infinite loop spaces, and define it inductively.

\begin{definition}
Define first an action of $\mathcal{U}^{\mathrm{op}}$ by
\[
d_{i}\ast1=\left\{
\begin{array}{ll}
1,&\text{if }i=0\\
0,&\text{otherwise}%
\end{array}
\right.
\]
\[
d_{i}\ast e_{j}e_{L}=\sum_{k}\left(  -1\right)  ^{i-k}\binom{i+\frac{j-i}{p}%
}{i-k}e_{i+\frac{j-i}{p}-\left(  p-1\right)  k}(d_{k}\ast e_{L}),
\leqno{\hbox{and}}
\]
where $e_{L}$ is any monomial in the $e$'s. 
\end{definition}

The reader may calculate that our Nishida action is encoded in the formal
power series identity
\[
d(u^{p-1})\ast\lbrack e(v)\cdot\underline{\quad}\,]=e(\varphi(v,u))[d(\varphi
(u,v)^{p-1})\ast\underline{\quad}\,],
\]
where $\varphi(a,b)$ is the function defined above. One can also check
straightforwardly that both $\widetilde{\mathcal{U}}$ and $\mathcal{U}$
inherit the action of $\mathcal{U}^{\mathrm{op}}$, considered either as subalgebras or
quotient algebras of $\widehat{\mathcal{U}}$.

\begin{theorem}
[proved in \fullref{sec-a}]\label{Nish-Adem}The formula for $\ast$
respects the Adem relations in $\smash{\mathcal{U}^{\mathrm{op}}}$ and hence defines a genuine
action of $\mathcal{K}^{\mathrm{op}}$ on $\smash{\widehat{\mathcal{U}}}$, $\smash{\widetilde
{\mathcal{U}}}$ and $\mathcal{U}.$ Furthermore, it also respects the Adem
relations in $\smash{\widehat{\mathcal{U}}}$, $\smash{\widetilde{\mathcal{U}}}$ and
$\mathcal{U}$, and hence induces an action of $\mathcal{K}^{\mathrm{op}}$ on
$\widehat{\mathcal{K}}$, 
$\widetilde{\mathcal{K}}$ and $\mathcal{K}.$
\end{theorem}

Note that our Nishida action is itself a map of coalgebras, as may easily be
verified by induction on length in $\widehat{\mathcal{U}}$ using the formal
power series formulation above of the Nishida action, since $e(t)$ (as noted
above) and $d(t)$ are both grouplike. This means that the induced
contragradient $\mathcal{A}$--action mentioned above turns the duals of
$\widehat{\mathcal{U}}$, $\widetilde{\mathcal{U}}$, $\mathcal{U}$,
$\widehat{\mathcal{K}}$, 
$\widetilde{\mathcal{K}}$, and $\mathcal{K}$ into unstable algebras over the
Steenrod algebra. Clearly from the theorem the natural maps between these are
maps over the Steenrod algebra.

We comment at this point on the relationship between the action we have
defined above, the traditional Nishida relations in topology, and the natural
Steenrod algebra action in invariant theory; the gist is that they are all in
agreement. Formulas compatible with those of our definition above, valid in
the homology of infinite loop spaces, can be derived from Theorem 9.4 (ii) of
May \cite{May}, so our formulas will agree with those of the traditional Nishida
relations when restricted to the Dyer--Lashof algebra, ie the homology of
$QS^{0}$. It is also not hard to check (as we did in the $2$--primary case in
\cite{PW1}) that $\mathcal{K}$ is isomorphic as a coalgebra to the mod $p$
Dyer--Lashof algebra modulo Bocksteins \cite{Kech1994}. Hence the dual to
$\mathcal{K}_{n,\ast}$ is isomorphic to the $n$--th Dickson algebra of $GL_{n}$
invariants \cite{Wilkerson}, and in fact this is an isomorphism of algebras
over the Steenrod algebra \cite[page 33f]{Wilkerson,CLM}. So we see that
our Nishida action will induce the same $\mathcal{A}$--action as that from
invariant theory. This could also be verified by direct calculations similar
to those in our proof of \fullref{slinv} below.

More broadly, we can now enlarge any inquiry comparing the duals of
$\widehat{\mathcal{U}}$, $\widetilde{\mathcal{U}}$, $\mathcal{U}$,
$\widetilde{\mathcal{K}}$, and $\mathcal{K}$ to algebras of invariants: Beyond
just comparing the dualization of their coalgebra structures to certain
algebra structures, we can actually compare unstable algebras over the
Steenrod algebra, since both they and algebras of invariants now have
independent defined structure as algebras over the Steenrod algebra. Despite
suggestions in the literature that such isomorphisms automatically respect the
$\mathcal{A}$--action, it seems to us that this still requires proof, and is
one of the most interesting features to ponder.

In the case of the projection of $\mathcal{U}$ onto $\mathcal{K}$, this has
already been explored by Kechagias \cite[Theorem 4.11]{Kech}, who proved that the
natural inclusion of the Dickson algebra of general linear group invariants
into the algebra of invariants under $T_{n}$, the upper triangular group with
arbitrary units on the diagonal, is dual to the coalgebra surjection of
$\mathcal{U}_{n,\ast}$ onto $\mathcal{K}_{n,\ast}$, as indicated in the
commutative diagram below. While the corresponding vertical maps are over
the Steenrod algebra by naturality, it is not yet clear to us whether the
isomorphism between the dual of $\mathcal{U}_{n,\ast}$ and the triangular
invariants respects their independent $\mathcal{A}$--actions.%
\[
\begin{CD}
S^{\tiny{T_{n}}} @<\cong<<  \mathcal{U}_{n}^{\ast}\\
@AAA                            @AAA\\
S^{\tiny{GL_{n}}} @<\cong<<  \mathcal{K}_{n}^{\ast}
\end{CD}
\]
As an example of our broader aims, we intend to enlarge the
situation above to a second commutative diagram:
\[
\begin{CD}
S^{\tiny{\widetilde{T}_{n}}}   @<\cong<\omega<  \widetilde{\mathcal{U}}_{n}^{\ast}\\
@A \tau AA                                              @AA \sigma A\\
S^{\tiny{\smash{\widetilde{SL}_{n}}}} @<\cong<<         \widetilde{\mathcal{K}}_{n}^{\ast}
\end{CD}
\]
Note that the entire first diagram maps by inclusions into the second, forming
a cube, and the reader may check that the faces of the cube joining the first
diagram to the second will also commute. Here $\widetilde{T}_{n}$ consists of
the upper triangular matrices whose diagonal elements are $\pm1$, and
$\smash{\widetilde{SL}_{n}}$ is the group of matrices whose determinants are $\pm1$. In this diagram, the horizontal isomorphisms are those provided by the
following two theorems. All the maps will arise from matching natural
generators for the algebras in question. \fullref{commute} below treats
the commutativity of this diagram.

\begin{theorem}
[proved in \fullref{sec-c}]\label{trianginv}The dual of $\widetilde
{\mathcal{U}}_{n,\ast}$ is isomorphic to the ring 
$\smash{\mathbb{F}_{p}[t_{1},\ldots
,t_{n}]^{\tiny{\widetilde{T}_{n}}}}$ of invariants under the action of the
group $\widetilde{T}_{n}$.
\end{theorem}

\begin{theorem}
[proved in \fullref{sec-c}]\label{slinv}The dual of $\widetilde
{\mathcal{K}}_{n,\ast}$ is isomorphic as an algebra over the Steenrod algebra
to $\smash{\mathbb{F}_{p}[t_{1},\ldots,t_{n}]^{\tiny{\widetilde{SL}_{n}}}}$, the ring of
invariants under the action of the group $\smash{\widetilde{SL}_{n}}$.
\end{theorem}

There are two main points to this theorem. The first is that we are dealing
with operations that lie below the radar of the classical lower-indexed
operations associated with the Steenrod and Dyer--Lashof algebras. The second
is that the Steenrod algebra structure is entirely determined by our Nishida
formula in $\widetilde{\mathcal{U}}.$ It is also the case that the Steenrod
algebra structure on the Dickson algebra $\mathbb{F}_{p}[t_{1},\ldots
,t_{n}]^{\tiny{GL_{n}}}$ is determined by the Nishida formula in $\mathcal{U}$. In work now in progress, we shall compute the relationships of other
quotient bi-algebras of $\widehat{\mathcal{U}}$ to algebras of invariants of
other subgroups of $GL_{n}$ that contain the upper triangular group
$\widehat{T}_{n}$ defined above.

We now turn to the commutativity of the diagram above. We begin by defining
the maps in the diagram. Identify $\smash{\widetilde{\mathcal{K}}_{n,\ast}}^{\ast
}=\mathbb{F}_{p}[t_{1},\ldots,t_{n}]^{\tiny{\smash{\widetilde{SL}_{n}}}}$ via the
isomorphism described in the proof of \fullref{slinv}. Define maps \begin{align*}
\sigma &  \co\smash{\smash{\widetilde{\mathcal{K}}_{n,\ast}}^{\ast}}\rightarrow\widetilde
{\mathcal{U}}_{n,\ast}^{\ast},\\
\tau &  \co\smash{\smash{\widetilde{\mathcal{K}}_{n,\ast}}^{\ast}}=\mathbb{F}_{p}[t_{1}%
,\ldots,t_{n}]^{\tiny{\widetilde{SL}_{n}}}\rightarrow\mathbb{F}_{p}[t_{1}%
,\ldots,t_{n}]^{\tiny{\widetilde{T}_{n}}},\text{ and}\\
\omega &  \co\widetilde{\mathcal{U}}_{n,\ast}^{\ast}\rightarrow\mathbb{F}%
_{p}[t_{1},\ldots,t_{n}]^{\tiny{\widetilde{T}_{n}}},\text{ \ as follows.}%
\end{align*}
The map $\sigma$ is the dual of the map that imposes Adem relations, the map
$\tau$ is the map induced by the inclusion $\widetilde{T}_{n}\subseteq
\smash{\widetilde{SL}_{n}}$, and the map $\omega$ is the isomorphism described in
the proof of \fullref{trianginv}.

\begin{theorem}
[proved in \fullref{sec-c}]\label{commute}We have $\omega\circ
\sigma=\tau.$
\end{theorem}

We conclude this section by listing two important self-maps and concomitant
properties of these algebras.

(1)\qua Let $\kappa$ denote multiplication by $e_{0}$. Since the element
$e_{0}=d_{0}$ satisfies $\Delta (  e_{0})  =e_{0}\otimes e_{0}$, its
dual, $\kappa^{\ast}$, is an algebra endomorphism on the duals of the various
algebras we have defined.

(2)\qua  The $p$--th power map on the dual algebras is known as the Frobenius map. Its dual $V$, known as the Verschiebung, is given on generators by $V\left(
e_{i}\right)  =e_{i/p}$, where we recall the convention that $e_{a}=0$ if $a$
is not a nonnegative integer. We extend the map multiplicatively to
products. Since $V(e(u))=e(u^{p})$, it is easy to check, using the power
series formulation of the Adem relations, that $V$ is well-defined on
$\widetilde{\mathcal{K}}$ and $\mathcal{K}$.

\begin{remark}
Relations between these maps and the Nishida action include:

(A)\qua $d_{i}\ast\kappa\left(  e_{J}\right)  =\kappa\left(  V\left(
d_{i}\right)  \ast e_{J}\right)  $, and

(B)\qua $d_{0}\ast e_{J}=V\left(  e_{J}\right)  .$

(Here $J=\left(  j_{1},\ldots,j_{n}\right)  $ is a multi-index and
$e_{J}=e_{j_{1}}\cdots e_{j_{n}}.$)
\end{remark}

\section[Proofs of Theorem 1.3, Theorem 1.9 and Theorem 1.18]{Proofs of \fullref{eademcoord}, \fullref{dfps} and
\fullref{Nish-Adem}}\label{sec-a}

\begin{proof}
[Proof of \fullref{eademcoord}]We will use the residue method of Jacobi
\cite{Jacobi} (see \cite[\S 1.1, 1.2]{Huang1,Huang2,Goulden}), and leave
straightforward calculations to the reader. If $f\left(  x\right)  =\sum
_{k}a_{k}x^{k}$ is a formal Laurent series (ie $k\in\mathbb{Z}$ is bounded
below) with coefficients in a ring with unity, define its \textquotedblleft
residue\textquotedblright\ $\operatorname*{res}_{x}\sum_{k}a_{k}x^{k}$ to be
$a_{-1}$. Then Jacobi's change of variables formula implies \cite[\S 1.1,
1.2]{Goulden} that if $y=g\left(  x\right)  $ is a formal power series with
coefficients in the same coefficient ring, and if $g\left(  0\right)  =0$ and
$g^{\prime}\left(  0\right)  $ is invertible in the coefficient ring, then
$\operatorname*{res}_{y}f\left(  y\right)  =\operatorname*{res}_{x}f\left(
g\left(  x\right)  \right)  g^{\prime}\left(  x\right)  $, where $g^{\prime}$
is the formal derivative.

First we will show that the full relations follow from the formal power series
identity form of May's relations.

In preparation, our setting is the ring $\big(  \widehat{\mathcal{U}}(
(  u) ) \big)  \big((v)\big)  $ of
formal Laurent series in $v$ with coefficients in the ring of formal Laurent
series in $u$ (with coefficients in $\smash{\widehat{\mathcal{U}}}$). Note that
$w=g\left(  v\right)  =v(v^{p-1}-u^{p-1})$ satisfies the hypotheses for a
Jacobi change of variables. Note too for use below that since $g$ is a very
simple multiplicatively invertible polynomial in $v$, one can compute its
powers in our ring, both positive and negative, by writing $g(  v)
^{m}=\big(  -u^{p-1}v(  1+(  -(  \frac{v}{u})
^{p-1}))  \big)  ^{m}$ and using the geometric/binomial power
series expansion for $(  1+x)  ^{m}$, for any $m\in\mathbb{Z}$.

Now using a Jacobi change of variables and the formal power series identity
for May's relations, we have, for any $i,j$ $\in\mathbb{Z}$,%
\begin{align*}
e_{i}e_{j}  &  =\operatorname*{res}\limits_{u}\left(  \operatorname*{res}%
_{w}\frac{e\left(  u\right)  e\left(  w\right)  }{u^{i+1}w^{j+1}}\right) \\
&  =\operatorname*{res}\limits_{u}\left(  \operatorname*{res}_{v}%
\frac{e\left(  u\right)  }{u^{i+1}}\frac{e\left(  v(v^{p-1}-u^{p-1})\right)
}{\left(  v(v^{p-1}-u^{p-1})\right)  ^{j+1}}\frac{d}{dv}\left(  v(v^{p-1}%
-u^{p-1})\right)  \right) \\
&  \sim\operatorname*{res}\limits_{u}\left(  \operatorname*{res}_{v}%
\frac{e(v)e(u(u^{p-1}-v^{p-1}))\left(  -u^{p-1}\right)  }{u^{i+1}\left(
v(v^{p-1}-u^{p-1})\right)  ^{j+1}}\right)  .
\end{align*}
From here a straightforward calculation of the latter using expanded formal
power series yields%
\[
\sum_{k}(-1)^{\frac{pk-i}{p-1}}\binom{k-j-1}{\frac{pk-i}{p-1}-j}%
e_{i+pj-pk}e_{k}\text{\ ,}%
\]
so the full relations follow from May's relations.

Now we turn to the converse, to prove the formal power series identity
relations under the hypothesis of the full relations. This time we prepare for
a Jacobi change of variables using the power series
\[
\rho\left(  w\right)  =-u\sum_{i\geq0}\left(  \frac{w}{u^{p}}\right)  ^{p^{i}%
}\text{,}%
\]
again with coefficients in $\widehat{\mathcal{U}}\left(  \left(  u\right)
\right)  $, and satisfying the Jacobi requirements since $\rho\left(
0\right)  =0$ and $\rho^{\prime}\left(  w\right)  =-u^{-\left(  p-1\right)  }$
(our characteristic is $p$), which is invertible in $\widehat{\mathcal{U}%
}\left(  \left(  u\right)  \right)  $. Note also that
\[
\left(  \frac{\rho\left(  w\right)  }{u}\right)  ^{p}=\frac{\rho\left(
w\right)  }{u}+\frac{w}{u^{p}},
\]
\[
w=\rho\left(  w\right)  ^{p}-u^{p-1}\rho\left(  w\right)  .
\leqno{\hbox{and thus}}
\]
In other words, $g\left(  \rho\left(  w\right)  \right)  =w$, so $\rho$ is the
composition inverse of $g$, and thus $\rho\left(  g\left(  v\right)  \right)
=v$ holds too \cite[\S 1.1]{Goulden}, which can easily be verified by direct
calculation and will also be needed below.

Now we assume the full relations as given, and will derive May's relations.
From the full relations and the calculations mentioned above, we have, for
every $i,j$ $\in\mathbb{Z}$,
\begin{align*}
e_{i}e_{j}  &  \sim\sum_{k}(-1)^{\frac{pk-i}{p-1}}\binom{k-j-1}{\frac
{pk-i}{p-1}-j}e_{i+pj-pk}e_{k}\text{\ }\\
&  =\operatorname*{res}\limits_{u}\left(  \operatorname*{res}_{t}%
\frac{e(t)e(u(u^{p-1}-t^{p-1}))\left(  -u^{p-1}\right)  }{u^{i+1}\left(
t(t^{p-1}-u^{p-1})\right)  ^{j+1}}\right)  .
\end{align*}
(Here we use $t$ to avoid confounding at this stage with the $v$ in our
desired final identity.) Next we make our change of variables $t=\rho\left(
w\right)  $, and use the polynomial equation above satisfied by $\rho\left(
w\right)  $, simplifying to
\[
e_{i}e_{j}\sim\operatorname*{res}\limits_{u}\Big(  \operatorname*{res}%
_{w}\frac{1}{u^{i+1}w^{j+1}}e(  \rho(  w)  )  e(
u(  u^{p-1}-\rho(  w)  ^{p-1})  )  \Big)  ,
\]
\[
e(  u)  e(  w)  \sim e(  \rho(  w)
)  e\big(  u(  u^{p-1}-\rho(  w)  ^{p-1}\big)
)  .
\leqno{\hbox{so}}
\]
Now we substitute $g(  v)  =v(v^{p-1}-u^{p-1})$ for $w$ in this
equality, and use the fact that $\rho(  g(  v)  )  =v$,
producing
\[
e(u)e\big(v(v^{p-1}-u^{p-1})\big)\sim e(v)e\big(u(u^{p-1}-v^{p-1})\big),
\]
as desired.
\end{proof}

\begin{proof}
[Proof of \fullref{dfps}]We begin with May's relations in $\widehat
{\mathcal{U}}$ in the form%
\[
e(u)e(v(v^{p-1}-u^{p-1}))\sim e(v)e(u(u^{p-1}-v^{p-1}))\text{.}%
\]
Now the map to the quotient $\mathcal{U}$ sends $e_{i\left(  p-1\right)  }$ to
$d_{i}$, and all other $e$'s to $0$, so the relations become
\[
d(  u^{p-1})  d\big(  v^{p-1}(  v^{p-1}-u^{p-1}\big)
^{p-1})  \sim d(  v^{p-1})  d\big(  u^{p-1}(
u^{p-1}-v^{p-1})  ^{p-1}\big)  ,
\]
and substituting $u$ for $u^{p-1}$ and $v$ for $v^{p-1}$ produces%
\[
d(u)d(v(v-u)^{p-1})\sim d(v)d(u(u-v)^{p-1}),
\]
as desired.
\end{proof}

In preparation for the proof of \fullref{Nish-Adem}, note the identity
\[
\varphi\left(  \varphi\left(  a,b\right)  ,\varphi\left(  c,b\right)  \right)
=\varphi\left(  \varphi\left(  a,c\right)  ,\varphi\left(  b,c\right)
\right)
\]
for the function $\varphi$ defined in the introduction.

\begin{proof}
[Proof of \fullref{Nish-Adem}]Let $\cdot$ denote the multiplication in
$\widehat{\mathcal{U}}$
and $\bullet$ be the multiplication in $\mathcal{U}^{\mathrm{op}}.$ For
the first part of the theorem, we must show that
\begin{align*}
&  \big(  d(\psi(v^{p-1},u^{p-1}))\bullet d(u^{p-1})\big)  \ast\big(
e(w)\cdot\underline{\quad}\,\big) \\
&  \sim\big(  d(\psi(u^{p-1},v^{p-1}))\bullet d(v^{p-1})\big)  \ast\big(
e(w)\cdot\underline{\quad}\,\big)  .
\end{align*}
We compute, assuming the result true inductively for lower length in
$\widehat{\mathcal{U}}$:%
\begin{align*}
&  \big(  d(\psi(v^{p-1},u^{p-1}))\bullet d(u^{p-1})\big)  \ast\big(
e(w)\cdot\underline{\quad}\,\big) \\
&  =d(\varphi(v,u)^{p-1})\ast\lbrack d(u^{p-1})\ast\big(  e(w)\cdot
\underline{\quad}\,\big)  ]\\
&  =d(\varphi(v,u)^{p-1})\ast\big[  e(\varphi(w,u))\cdot\big(  d\big(
\varphi(u,w)^{p-1}\big)  \ast\underline{\quad}\,\big)  \big] \\
&  =e\big(  \varphi(\varphi(w,u),\varphi(v,u))\big)  \cdot\big[  d\big(
\varphi(\varphi(v,u),\varphi(w,u))^{p-1}\big)  \ast\big(  d\big(
\varphi(u,w)^{p-1}\big)  \ast\underline{\quad}\,\big)  \big] \\
&  =e\big(  \varphi(\varphi(w,u),\varphi(v,u))\big)  \cdot\big[  d\big(
\varphi(\varphi(v,w),\varphi(u,w))^{p-1}\big)  \ast\big(  d\big(
\varphi(u,w)^{p-1}\big)  \ast\underline{\quad}\,\big)  \big] \\
&  =e\big(  \varphi(\varphi(w,u),\varphi(v,u))\big)  \cdot\Big[\!  \Big(
d\big(  \varphi(\varphi(v,w),\varphi(u,w))^{p-1}\big)  \bullet d\big(
\varphi(u,w)^{p-1}\big)  \Big)  \ast\underline{\quad}\,\Big] \\
&  =\!\smash{e\big(  \varphi(\varphi(w,u),\varphi(v,u))\big)  \cdot\Big[ \! \Big(
\!d\big(  \psi(\varphi(v,w)^{p-1},\varphi(u,w)^{p-1})\big)  \bullet d\big(
\varphi(u,w)^{p-1}\big)\!  \Big)  \ast\underline{\quad}\Big]} \\
&  \sim\! e\big(  \varphi(\varphi(w,v),\varphi(u,v))\big)  \cdot\Big[\!
\Big(  d\big(  \psi(\varphi(u,w)^{p-1},\varphi(v,w)^{p-1})\big)  \bullet
d\big(  \varphi(v,w)^{p-1}\big) \! \Big)  \!\ast\underline{\quad}\,\Big].
\end{align*}
We may now reverse these steps to see that this last term is equal to the desired
$$\big(
d(\psi(u^{p-1},v^{p-1}))\bullet d(v^{p-1})\big)  \ast\big(  e(w)\cdot
\underline{\quad}\,\big).$$
For the second part of \fullref{Nish-Adem}, we must show that%
\[
d(u^{p-1})\ast\big[  e(v)\cdot e(\varphi(w,v))\cdot\underline{\quad}\,\big]
\sim d(u^{p-1})\ast\big[  e(w)\cdot e(\varphi(v,w))\cdot\underline{\quad
}\,\big]  .
\]
We compute%
\begin{align*}
&  d(u^{p-1})\ast\big[  e(v)\cdot e(\varphi(w,v))\cdot\underline{\quad
}\,\big] \\
&  =e(\varphi(v,u))\cdot\big[  d\big(  \varphi(u,v)^{p-1}\big)
\ast\big(  e\big(  \varphi(w,v)\big)  \cdot\underline{\quad}\,\big)
\big] \\
&  =e(\varphi(v,u))\cdot\Big[  e\big(\varphi\big(\varphi(w,v),\varphi(u,v)\big)\cdot
\Big(  d\big(\varphi\big(\varphi(u,v),\varphi(w,v)\big)\big)  ^{p-1}\ast\underline
{\quad}\,\Big)\Big] \\
&  =\smash{\big[  e(\varphi(v,u))\cdot e\big(\varphi(\varphi(w,v),\varphi(u,v)\big)\big]
\cdot\Big(  d\big(\varphi\big(\varphi(u,v),\varphi(w,v)\big)\big)  ^{p-1}\ast
\underline{\quad}\,\Big)}\\
&  =\big[  e(\varphi(v,u))\cdot e\big(\varphi(\varphi(w,u),\varphi(v,u)\big)\big]
\cdot\Big(  d\big(\varphi(\varphi(u,w),\varphi(v,w)\big)\big)  ^{p-1}\ast
\underline{\quad}\,\Big)\\
&  \sim\big[  e(\varphi(w,u))\cdot e\big(\varphi(\varphi(v,u),\varphi
(w,u)\big)\big]  \cdot\Big(  d\big(\varphi(\varphi(u,w),\varphi(v,w)\big)\big)
^{p-1}\ast\underline{\quad}\,\Big).
\end{align*}
By reversing these steps, this last term is equal to
$d(u^{p-1})\ast\big[  e(w)\cdot e(\varphi(v,w))\cdot\underline{\quad}\,\big]  $,
as desired.
\end{proof}

\section[Proofs of Proposition 1.12, Theorem 1.13 and Theorem 1.14]{\label{sec-b}Proofs of \fullref{noncongruenttrivialproducts},
\fullref{basis} and \fullref{negredundant}}

To begin this section, we collect some facts about the full relations in
$\widehat{\mathcal{U}}_{2,\ast}$. We use the notation $\theta(i,j)$ from the
introduction for the right side of any full relation.

\begin{lemma}
\label{Adem-index-congruence}For any $i,j\in\mathbb{Z}$, if $e_{l}e_{k}$
appears with a nonzero coefficient in $\theta(i,j)$, then $k\equiv i$
$\operatorname{mod}\left(  p-1\right)  $ and $l\equiv j$ $\operatorname{mod}%
\left(  p-1\right)  $.
\end{lemma}

(So the second index in each term that appears on the right-hand side
of any full relation is congruent mod $\left(  p-1\right)  $ to the first
index on the left-hand side, and similarly for the other pair.)

\begin{proof}
The fraction that occurs as the exponent of $-1$ must be an integer, so
$k\equiv i$ $\operatorname{mod}\left(  p-1\right)  $. The second congruence
then follows from $i+pj=l+pk$.
\end{proof}

Thus we have the following lemma.

\begin{lemma}
\label{Adem-congruence-classes}For any $i,j\in\mathbb{Z}$, if $e_{l}e_{k}$
appears with a nonzero coefficient in $\theta(i,j)$, then $l-k\equiv-\left(
i-j\right)  $ $\operatorname{mod}\left(  p-1\right)  $.
\end{lemma}

(So the difference of indices in each term that appears on the
right-hand side of any full relation is congruent mod $\left(  p-1\right)  $
to the negation of the corresponding index difference on the left. Thus
$\widehat{\mathcal{U}}_{2,\ast}$ splits into two subspaces, generated
respectively by those $e_{i}e_{j}$ for which either $i\equiv j$ or
$i\not \equiv j$ $\operatorname{mod}\left(  p-1\right)  $, and every full
relation involves terms that lie in only one of these subspaces.)

\begin{lemma}
\label{Adem-inadmissibles}For any $i,j\in\mathbb{Z}$, if $i>j$ and $e_{l}%
e_{k}$ appears with a nonzero coefficient in $\theta(i,j)$, then $l\leq k.$
\end{lemma}

(That is, any full relation for an inadmissible monomial rewrites the
monomial in terms of admissibles.)

\begin{proof}
For the denominator of the binomial coefficient to be nonnegative, we must
have $pk\geq i+(p-1)j$, so $i+pj-pk\leq j$, whence $k\geq
i+pj-pk+(k-j)=l+(k-j).$ Since $i>j$, the first inequality also yields
$pk>pj$, whence $k>j.$ Combining these facts, we see that $l\leq k.$
\end{proof}

\begin{lemma}
\label{Adem-numerators}For any $i,j\in\mathbb{Z}$, if $i\leq j$ and
$e_{l}e_{k}$ appears with a nonzero coefficient in $\theta(i,j)$, then the
numerator of the binomial coefficient is negative. Conversely, if $i>j$ and
$e_{l}e_{k}$ appears with a nonzero coefficient, then the numerator of the
binomial coefficient is nonnegative.
\end{lemma}

(That is, any full relation for an admissible monomial produces
exclusively terms on the right with negative numerators in their binomial
coefficients, and for an inadmissible monomial produces exclusively terms with
nonnegative numerators.)

\begin{proof}
Consider first $i\leq j$. If the numerator were nonnegative, $k\geq j+1$, and
the binomial coefficient is nonzero, then $k-j-1\geq\frac{pk-i}{p-1}-j$, from
which we obtain $i-(p-1)\geq k.$ Then since $i\leq j$, we have $j-(p-1)\geq
k\geq j+1$, a contradiction. Hence $k<j+1$, ie the numerator is
negative$.$

Next consider $i>j$, and the binomial coefficient nonzero, so the denominator
is nonnegative, ie $pk-i\geq pj-j$. Then $k-j\geq\upnfrac{i-j}{p}>0$, so
$k-j-1\geq0$, ie the numerator is nonnegative.
\end{proof}

\begin{lemma}
\label{Adem-edge}Let $1\leq b\leq p-1.$ Then 
$\theta(i,i-b)=\theta (  e_{i}e_{i-b})  =0$ for 
any $i\geq b$.
\end{lemma}

(This expresses an important \textquotedblleft edge
effect\textquotedblright, in which Adem relations for \textquotedblleft nearly
admissible\textquotedblright\ inadmissibles are zero on the right side.)

\begin{proof}
We compute
\[
e_{i}e_{i-b}\sim\sum_{k}(-1)^{\frac{pk-i}{p-1}}\binom{k-i+b-1}{\frac
{pk-i}{p-1}-i+b}e_{i+p(i-b)-pk}e_{k}.
\]
Since $b\geq1$, the left side is inadmissible, so by \fullref{Adem-numerators} the numerator of the binomial coefficient is
nonnegative. Thus for it to be nonzero, we must have $k-i+b-1$ $\geq
\frac{pk-i}{p-1}-i+b$, whence $i-(p-1)\geq k.$ But also $k\geq i-b+1.$ Since $i-(p-1)-\left(  i-b+1\right)  =b-p$, there are no values of $k$ that
satisfy both inequalities if $1\leq b\leq p-1.$
\end{proof}

In preparation for the proof of \fullref{basis}, we follow Le Minh Ha
\cite{Le}, who proved algebraically that the Adem relations for admissibles in
the Steenrod algebra are redundant, using the self-map of the Steenrod algebra
dual to multiplication by $\xi_{1}$ in the dual Steenrod algebra. We imitate
his self-map in our setting, defining
\begin{gather*}
\eta\co\widehat{\mathcal{U}}_{2,\ast}\longrightarrow\widehat{\mathcal{U}%
}_{2,\ast+p-1}\\
\tag*{\hbox{by}}
\eta(e_{i}e_{j})=e_{i+p-1}e_{j}+e_{i-(p-1)^{2}}e_{j+p-1}\text{ for }i,j\geq0.
\end{gather*}
NB: Like the formula for $\theta(  e_{i}e_{j})  $, this formula
only applies when $i,j\geq0$. We need to be extremely careful to pay attention
to this, and will say \textquotedblleft $\eta$ is defined by
formula\textquotedblright\ for emphasis when this is the case.

\begin{lemma}
\label{eta-Adem}Consider $e_{i}e_{j}$ for $i\geq (  p-1 )  ^{2}$ and
$j\geq0.$

{\rm (a)}\qua If $0\leq i\leq j$ (ie $e_{i}e_{j}$ is admissible), then $\theta (
\eta (  e_{i}e_{j} )   )  =\eta (  \theta (  e_{i}%
e_{j} )   )  .$

{\rm (b)} If $0\leq j<i$ (ie $e_{i}e_{j}$ is inadmissible), and $i-j\leq\frac
{i}{p}+p-1$, then $\theta (  \eta (  e_{i}e_{j} )   )
=\eta (  \theta (  e_{i}e_{j} )   )  .$
\end{lemma}

\begin{proof}
First note, regarding the left side of the equality, that since $i\geq (
p-1)  ^{2}$, both terms in the formula above for $\eta(e_{i}e_{j})$ have
entirely nonnegative indices, and thus $\theta$ is defined on them by the
formula for Adem relations in its definition.

Second, regarding the right side of the equality, in both cases (a) and (b) we
shall check that if $e_{l}e_{k}$ is a term in $\theta (  e_{i}%
e_{j})  $ with nonzero binomial coefficient, then $l$ and $k$ are
nonnegative, and hence $\eta$ is defined by formula on $e_{l}e_{k}$. Indeed,
in case (a), by \fullref{Adem-numerators}, $k\leq j$. So $l=i+pj-pk\geq
i\geq0$. And $\frac{pk-i}{p-1}-j\geq0$, whence $pk\geq i+\left(  p-1\right)
j$, so $k\geq0.$ In case (b), by \fullref{Adem-numerators}, $k\geq j+1\geq
0$. So also $\frac{pk-i}{p-1}-j\leq k-j-1$, thus $pk-i\leq\left(  p-1\right)
k-\left(  p-1\right)  $, and so $k\leq i-\left(  p-1\right)  $. Then
\begin{align*}
l=i+pj-pk  &  =\left(  p+1\right)  i-p\left(  i-j\right)  -pk\\
&  \geq\left(  p+1\right)  i-p\left(  i-j\right)  -pi+p\left(  p-1\right) \\
&  =i-p\left(  i-j\right)  +p\left(  p-1\right) \\
&  \geq i-\left(  i+\left(  p-1\right)  p\right)  +p\left(  p-1\right)  =0,
\end{align*}
using the additional hypothesis for the last inequality.

Now, since we have shown that all terms on both sides of the claimed equality
$$\theta (\eta(e_{i}e_{j}))=\eta(\theta(e_{i}e_{j}))$$
will be calculated using the
defining formulas for $\theta$ and $\eta$ on basis elements, it remains to
show that these match up. We leave the details to the reader, noting that it
boils down first to combining four terms to three via Pascal's ordinary
binomial coefficient identity, followed by the mod $p$ identity $\binom{M}%
{N}\equiv\binom{M-p}{N}+\binom{M-p}{N-p}$ also used by Le Minh Ha \cite{Le}.
\end{proof}

\begin{lemma}
\label{lesser-excess}For $i,j\geq0$, if $i\leq j$, then in the quotient
$\widehat{\mathcal{K}}$, the admissible $e_{i}e_{j}$ can be rewritten by
application of Adem relations as
\[
e_{i}e_{j}\text{ }=\gamma e_{i}e_{j}+\text{a sum of admissibles of lesser
excess,}%
\]
with coefficient $\gamma=0$ if $i\not \equiv j$ $\operatorname{mod}\left(
p-1\right)  $, and $\gamma=\left(  -1\right)  ^{i}=\left(  -1\right)  ^{j}$ if
$i\equiv j$ $\operatorname{mod}\left(  p-1\right)  $.
\end{lemma}

\begin{proof}
In the quotient $\widehat{\mathcal{K}}$ by the Adem relations, we have
\[
e_{i}e_{j}=\sum_{k}(-1)^{\frac{pk-i}{p-1}}\binom{k-j-1}{\frac{pk-i}{p-1}%
-j}e_{l}e_{k},
\quad \text{where }%
l=i+pj-pk.
\]
Since $i\leq j$, by \fullref{Adem-numerators}, for nonzero terms we have
$k\leq j$. By \fullref{Adem-index-congruence}, there is no term with $k=j$
unless $j=k\equiv i$ $\operatorname{mod}\left(  p-1\right)  $), in which case
the coefficient is $\left(  -1\right)  ^{-j}$, the sign claimed for $\gamma$.

Now the remaining terms in the sum are either admissible of lesser excess,
since $k<j$ and $i+pj=l+pk$, or inadmissible. If inadmissible, ie $l>k$, we
proceed as follows for any nonzero term.

First, since we must have
\begin{align*}
&\quad\frac{pk-i}{p-1}-j\geq0,
\\
\tag*{\hbox{we get}}
pk&\geq i+(p-1)j.\\
\tag*{\hbox{Then}}
l  &  =i+pj-pk\\
&  \leq i+pj-\left(  i+(p-1)j\right)  =j.
\end{align*}
Now applying another Adem relation, from \fullref{Adem-inadmissibles} we
have
\[
e_{l}e_{k}=\sum_{m}(-1)^{\frac{pm-l}{p-1}}\binom{m-k-1}{\frac{pm-l}{p-1}%
-k}e_{i+pj-pm}e_{m},
\]
where the nonzero terms on the right are all admissible. By \fullref{Adem-numerators}, $m>k.$ Also, $m-k-1\geq\smash{\frac{pm-l}{p-1}}-k,$ whence
$
l-(p-1)\geq m,
$
so%
$
j-(p-1)\geq m\text{ from above,}%
$
and thus
$
j>m.
$
Write $n=i+pj-pm$, so we now have $n>i$. Thus $m-n<j-i$, so $e_{n}e_{m}$ is an
admissible of lesser excess than $e_{i}e_{j}$.
\end{proof}

\begin{proof}
[Proof of \fullref{noncongruenttrivialproducts}]Clearly we may assume
$i,j\geq0$. We begin with part (a). First, consider $i>j$. \fullref{Adem-congruence-classes} and \fullref{Adem-inadmissibles}, along with the
hypothesis $i$ $\not \equiv $ $j$ $\operatorname{mod}\left(  p-1\right)  $,
ensure that in the quotient $\widehat{\mathcal{K}}$ by the Adem relations,
$e_{i}e_{j}$ is a linear combination of terms $e_{l}e_{k}$ with $l$
$\not \equiv $ $k$ $\operatorname{mod}\left(  p-1\right)  $ and $l\leq k.$
This reduces the proof to considering terms for which $i\leq$ $j.$ So let
$0\leq i\leq$ $j$ and $i$ $\not \equiv $ $j$ $\operatorname{mod}\left(
p-1\right)  .$ Then from \fullref{Adem-congruence-classes} and
\fullref{lesser-excess}, $e_{i}e_{j}$ can be written as a sum of admissibles of
lower excess satisfying the same noncongruence condition. By Fermat's method
of descent the proof is complete.

For part (b), we first note that thanks to part (a) we need only consider the
case where $i\equiv j\operatorname{mod}\left(  p-1\right)  $. Then 
\fullref{Adem-index-congruence}, \fullref{Adem-congruence-classes}, and
\fullref{Adem-inadmissibles} ensure that in the quotient $\widehat{\mathcal{K}}$
by the Adem relations, $e_{i}e_{j}$ is a linear combination of terms
$e_{l}e_{k}$ with $l$ $\equiv$ $k$ $\operatorname{mod}\left(  p-1\right)  $,
$0\leq l\leq k$, and $l,k$ odd. Now \fullref{lesser-excess} writes
$2e_{l}e_{k}$ as a sum of terms of lesser excess, that by our Lemmas also
satisfy the same conditions of index congruence, admissibility, and oddness.
Again by descent we are finished.
\end{proof}

\begin{proof}
[Proof of \fullref{basis}]Since we are working to analyze $\widetilde
{\mathcal{K}}$, a quotient of $\widetilde{\mathcal{U}}$, henceforth all
subscripts on $e$'s shall be even. According to \fullref{Adem-inadmissibles}
and \fullref{noncongruenttrivialproducts}, we need only prove that the
elements $e_{i_{1}}\cdots e_{i_{n}\text{ }}$for which all $i_{k\text{ }}$ are
even, $0\leq i_{1}\leq\cdots\leq i_{n}$, and $i_{k\text{ }}\equiv
i_{k-1\text{ }}$ $\operatorname{mod}\left(  p-1\right)  $ for $2\leq k\leq n$,
are linearly independent in $\widetilde{\mathcal{K}}$. This proof will use
the method of Le Minh Ha \cite{Le}. From a little linear algebra and \fullref{Adem-congruence-classes} and \fullref{lesser-excess}, it is enough to
show that in $\widetilde{\mathcal{U}}$, the Adem relations for $e_{i}e_{j}$,
for $0\leq i\leq j$ and $i-j\equiv0$ $\operatorname{mod}\left(  p-1\right)
$, are consequences of the relations for which $i>j.$ Let $\Theta\left(
i,j\right)  =e_{i}e_{j}-\theta (  e_{i}e_{j})  $ for $i,j\geq0$. Our goal is thus to prove that in $\widetilde{\mathcal{U}}$ every
$\Theta\left(  i,j\right)  $ for which $j-i\equiv0$ $\operatorname{mod}%
\left(  p-1\right)  $ and $i\leq j$ ($i,j$ even) is a linear combination of
 $\Theta\left(  l,k\right)  $'s for which $l>k.$ Thus we consider
$\Theta\left(  i,i+s\right)  $ for $i,s\geq0$, and proceed by induction on
$s$. Since $s\equiv0$ $\operatorname{mod}\left(  p-1\right)  $, we may write
$s=c\left(  p-1\right)  $ and induct on $c\geq0$.

{\bf Case One (base step)}\qua Consider $s=c(p-1)$, where $0\leq c\leq p.$ We have
\[
\theta (  e_{i}e_{j})  =\theta (  e_{i}e_{i+s})  =\sum
_{k}(-1)^{\frac{pk-i}{p-1}}\binom{k-i-s-1}{\frac{pk-i}{p-1}-i-s}%
e_{(p+1)i+ps-pk}e_{k}.
\]
For the binomial coefficient to be nonzero, we must have (using \fullref{Adem-numerators} to know that the numerator is negative)
\begin{align*}
i+s  &  \geq k\geq i+s-\frac{s}{p}\\
&  \geq i+c(p-1)-\frac{c(p-1)}{p}\\
&  =i+c(p-1)-c+\frac{c}{p}.\\
\tag*{\hbox{So}}
i+c(p-1)&\geq k\geq i+c(p-1)-c+\frac{c}{p}.
\end{align*}
By \fullref{Adem-index-congruence}, there is one possible value for $k$
 in this range, namely $k=i+c(p-1)$, unless $c=p$, in which case $k=i+\left(
c-1\right)  (p-1)$ is also a possibility. For $k=i+c(p-1)$, we get the term
\begin{align*}
&  (-1)^{\frac{pk-i}{p-1}}\binom{k-i-s-1}{\frac{pk-i}{p-1}-i-s}%
e_{(p+1)i+ps-pk}e_{k}\\ &  =(-1)^{i+pc}\binom{-1}{c}e_{i}e_{i+c(p-1)}\\
&  =(-1)^{i}e_{i}e_{i+s}.
\end{align*}
Thus since $i$ is even, we have shown that for $0\leq c<p$, the Adem relation
for $e_{i}e_{i+c\left(  p-1\right)  }$ reduces to
\[
\Theta\left(  i,i+s\right)  =e_{i}e_{i+s}-e_{i}e_{i+s}=0,
\]
while for $c=p$ there is one possible remaining surviving term, with
$k=i+(p-1)^{2}$, a multiple of $\smash{e_{i+p^{2}-p}e_{i+(p-1)^{2}}}$. But
$\smash{\theta (  e_{i+p^{2}-p}e_{i+(p-1)^{2}})}  =0$ by \fullref{Adem-edge}, so $\smash{e_{i+p^{2}-p}e_{i+(p-1)^{2}}=\Theta (  i+p^{2}%
-p,i+(p-1)^{2})}  $, expressing the dependence of\break $\Theta\left(
i,i+s\right)  $ on an inadmissible relation, as desired.

{\bf Case Two (inductive step)}\qua Let $c\geq p+1.$ Inductively, we assume that
for $0\leq q<c$ we have
\[
\Theta(a,a+q(p-1))=\sum\limits_{m\in M_{q}}\gamma_{m}\Theta
(a+mp(p-1),a+(q-m)(p-1)),
\]
where $\gamma_{m}$ is a scalar and
\[
M_{q}=\{m|mp\leq q\leq m(p+1)-1\}.
\]
We note that the relations appearing on the right hand side of this formula
are inadmissible, and that the Case One results do satisfy this assumption.

By \fullref{eta-Adem}(a), we have
\begin{multline*}
\Theta(i,i+c(p-1))\\
=\eta\Theta(i+(p-1)^{2},i+(c-1)(p-1))-\Theta
(i+p(p-1),i+(c-1)(p-1)).
\end{multline*}
We focus first on the second term on the right side. Since $c-1\geq p$, by the
inductive hypothesis we have
\begin{multline*}
\Theta(i+p(p-1),i+(c-1)(p-1))\\
=\sum\limits_{m^{\prime}\in M_{q^{\prime}%
}^{\prime}}\gamma_{m^{\prime}}\Theta(a^{\prime}+m^{\prime}p(p-1),a^{\prime
}+(q^{\prime}-m^{\prime})(p-1)),
\end{multline*}
where $a^{\prime}=i+p(p-1),$ $q^{\prime}=c-1-p,$ and
\[
M_{q^{\prime}}^{\prime}=\{m^{\prime}|m^{\prime}p\leq q^{\prime}\leq m^{\prime
}(p+1)-1\}.
\]
\begin{multline*}
\tag*{\hbox{Here}}\Theta(a^{\prime}+m^{\prime}p(p-1),a^{\prime}+(q^{\prime}-m^{\prime
})(p-1))\\
=\Theta(i+(m^{\prime}+1)p(p-1),i+(c-1-m^{\prime})(p-1)).
\end{multline*}
And $m^{\prime}p\leq c-1-p\leq m^{\prime}(p+1)-1,$ so $(m^{\prime}+1)p+1\leq
c\leq(m^{\prime}+1)(p+1)-1.$

Now focusing on the first term on the right side, we can again use the
inductive hypothesis to obtain
\begin{multline*}
\eta\Theta(i+(p-1)^{2},i+(c-1)(p-1))\\
=\eta\sum\limits_{m^{\prime\prime}\in
M_{q^{\prime\prime}}^{\prime\prime}}\gamma_{m^{\prime\prime}}\Theta
(a^{\prime\prime}+m^{\prime\prime}p(p-1),a^{\prime\prime}+(q^{\prime\prime
}-m^{\prime\prime})(p-1)),
\end{multline*}
where $a^{\prime\prime}=i+(p-1)^{2},$ $q^{\prime\prime}=c-p,$ and
\[
M_{q^{\prime\prime}}^{\prime\prime}=\{m^{\prime\prime}\,|\,m^{\prime\prime}p\leq
q^{\prime\prime}\leq m^{\prime\prime}(p+1)-1\}.
\]
One checks that $(a^{\prime\prime}+m^{\prime\prime}p(p-1),a^{\prime\prime
}+(q^{\prime\prime}-m^{\prime\prime})(p-1))$ satisfies the hypotheses of
\fullref{eta-Adem}(b). Hence
\begin{multline*}
  \eta\Theta(a^{\prime\prime}+m^{\prime\prime}p(p-1),a^{\prime\prime
}+(q^{\prime\prime}-m^{\prime\prime})(p-1))\\
  =\Theta(i+(m^{\prime\prime}+1)p(p-1),i+(c-1-m^{\prime\prime})(p-1))\\
  +\Theta(i+m^{\prime\prime}p(p-1),i+(c-m^{\prime\prime})(p-1)).
\end{multline*}
Considering the first term on the right here, we have $m^{\prime\prime}p\leq
c-p\leq m^{\prime\prime}(p+1)-1,$ so $(m^{\prime\prime}+1)p\leq c\leq
(m^{\prime\prime}+1)(p+1)-2$, and therefore the term is inadmissible and meets
the inductive requirement. For the second term, $m^{\prime\prime}%
p<(m^{\prime\prime}+1)p\leq c;$ and either $c\leq m^{\prime\prime}(p+1)-1$
or $m^{\prime\prime}(p+1)\leq c\leq(m^{\prime\prime}+1)(p+1)-2.$ In the first
instance, the term is inadmissible and the inductive assumption is met. In the
latter case, we have an admissible for which $0\leq c-m^{\prime\prime
}(p+1)\leq p-1,$ so by Case 1, $\Theta(i+m^{\prime\prime}p(p-1),i+(c-m^{\prime
\prime})(p-1))=0.$

We conclude that the elements of the form $e_{a}e_{a+q}$ for $a$ even,
$q\geq0$, and $q$ divisible by $ p-1$, are linearly independent in
$\widetilde{\mathcal{K}}_{2,\ast}.$ Both parts of the theorem follow.
\end{proof}

\begin{proof}
[Proof of \fullref{negredundant}]We will use the notation $\Theta(i,j)$ of
the previous proof for the Adem relations, but extended to the full relations
for $i,j\in\mathbb{Z}$ via $\Theta(i,j)=e_{i}e_{j}-\theta(i,j)$, where
$\theta(i,j)$ is our notation from the introduction for the right side of any
full relation (and differs in general, we recall, from $\theta\left(
e_{i}e_{j}\right)  $, which is zero when either $i$ or $j$ is negative). Now
for any $i,j\in\mathbb{Z}$, we have $\wtilde{\alpha}^{n}\left(
\Theta(i,j)\right)  =\Theta\left(  i-2n,j-2n\right)  $. This is essentially
contained in the fact that the formal power series formulation of the Adem
relations clearly commutes with $\wtilde{\alpha}$ since $\widetilde{\alpha
}\left(  e(u)\right)  =u^{2}e(u)$, but it may be verified directly by
comparing the explicit formulas.

We will first prove our claim for the negative Adem relations on
$\widetilde{\mathcal{U}}$, ie we will show that for any relation
$\Theta(i,j)$ on $\widetilde{\mathcal{U}}$ (ie $i,j$ even) for which $i<0$
or $j<0$, the relation is already trivial once the Adem relations on
$\widetilde{\mathcal{U}}$ are imposed, ie it is trivial in $\widetilde
{\mathcal{K}}$. First, since any nonzero term in the formula for $\theta(i,j)$
must satisfy $i+\left(  p-1\right)  j\leq pk\leq i+pj$, we see that if $j<0$,
there are none, so $\Theta(i,j)$ is trivial in $\widetilde{\mathcal{K}}$.
Furthermore, if $i\not \equiv j$ $\operatorname{mod}\left(  p-1\right)  $,
then all terms in $\theta(i,j)$ also have their indices noncongruent
$\operatorname{mod}\left(  p-1\right)  $ from \fullref{Adem-congruence-classes}, so by \fullref{noncongruenttrivialproducts} all such terms are also zero in
$\widetilde{\mathcal{K}}$. Thus it remains only to consider the case when
$i<0$, $j\geq0$. and $i\equiv j$ $\operatorname{mod}\left(  p-1\right)  $.

In this case, we first examine $\Theta(-i,j-2i),$which is also a relation on
$\widetilde{\mathcal{U}}$ satisfying the same congruence condition on its
entries, but now with both entries positive. Thus from the proof of \fullref{basis}, if it is not already inadmissible it can be written as a linear
combination of inadmissibles, ie in $\widetilde{\mathcal{U}}$,
\[
\Theta(-i,j-2i)=\sum_{k>l}\gamma_{k,l}\Theta\left(  k,l\right)  .
\]
Now we apply $\wtilde{\alpha}^{-i}$ (recall that $-i>0$), which from our
observation above, about how $\wtilde{\alpha}$ passes through $\Theta$,
yields
\[
\Theta(i,j)=\sum_{k>l}\gamma_{k,l}\Theta\left(  k+2i,l+2i\right)
\]
in $\widetilde{\mathcal{U}}$. Now if $l+2i<0$, we know from above that
$\Theta\left(  k+2i,l+2i\right)  =0$ in $\widetilde{\mathcal{U}}$, and if
$l+2i\geq0$, then since $k>l$, we see that $\Theta\left(  k+2i,l+2i\right)  $
is one of the Adem relations already imposed in the definition of
$\widetilde{\mathcal{K}}$.

To extend our claim to cover the negative relations on $\widehat{\mathcal{U}}%
$, notice that if either $i$ or $j$ is odd, then from \fullref{Adem-index-congruence} and \fullref{Adem-congruence-classes}, any nonzero
term in the formula for $\theta(i,j)$ must have one of its indices odd as
well, so from \fullref{noncongruenttrivialproducts} all such terms are
zero once the Adem relations are imposed; thus the relation $\Theta(i,j)$ is
trivial after imposition of the Adem relations on $\widehat{\mathcal{U}}$.
\end{proof}

\section[Proofs of Theorem 1.19, Theorem 1.20 and Theorem 1.21]{\label{sec-c}Proofs of \fullref{trianginv}, \fullref{slinv} and
\fullref{commute}}

\begin{proof}
[Proof of \fullref{trianginv}]We have seen that the coalgebra primitives
in $\widetilde{\mathcal{U}}_{n,\ast}$ are the monomial basis elements
$e_{0}^{n}$ and $\smash{e_{0}^{a}e_{2}e_{0}^{n-a-1}}$, for $0\leq a\leq n-1.$ Let
 $\wtilde{v}_{n,a}=$ $(  \smash{e_{0}^{a}e_{2}e_{0}^{n-a-1})  ^{\ast
}}$, the dual element to $\smash{e_{0}^{a}e_{2}e_{0}^{n-a-1}}$. Note that the
topological degree of $\wtilde{v}_{n,a}$ is $4p^{a}.$ Also note
that since $\Delta (  e_{0}^{n})  =e_{0}^{n}\otimes e_{0}^{n}$,
 its dual is $1$, the algebra identity in $\smash{\widetilde{\mathcal{U}%
}_{n,\ast}^{\ast}}.$ The multiplication in $\smash{\widetilde{\mathcal{U}}%
_{n,\ast}^{\ast}}$ is commutative and obeys the usual degree convention for
products. On basis elements of $\smash{\widetilde{\mathcal{U}}_{n,\ast}}$, the
correspondence%
\[
e_{2i_{1}}\cdots e_{2i_{n}}\mapsto\wtilde{v}_{n,0}^{i_{1}}\cdots
\wtilde{v}_{n,n-1}^{i_{n}}%
\]
provides a bijection of graded vector spaces from $\widetilde{\mathcal{U}%
}_{n,\ast}$ to $\mathbb{F}_{p}[  \wtilde{v}_{n,0},\ldots
,\wtilde{v}_{n,n-1}]  .$ It follows from calculations of M\`ui
\cite{Mui} and Kechagias \cite[Corollary 4.22]{Kech1993,Kech1994} that the equality $\smash{\mathbb{F}_{p}
[t_{1},\ldots,t_{n}]^{\tiny{\widetilde{T}_{n}}}=\mathbb{F}_{p}[  \widetilde
{V}_{1},\ldots,\widetilde{V}_{n}]}  $ holds, where the degree of
 $\widetilde{V}_{i}$ (the square of M\`ui's invariant $V_{i}$) is
 $4p^{i-1}$, so we obtain our desired result by mapping one set of generators
to the other.
\end{proof}

\begin{proof}
[Proof of \fullref{slinv}]First we calculate the structure of the dual of
$\smash{\widetilde{\mathcal{K}}_{n,\ast}}$ as an algebra over the Steenrod algebra.
The coalgebra primitives in $\smash{\widetilde{\mathcal{K}}_{n,\ast}}$ $ $are
$e_{2}^{n}$, and $e_{0}^{a}e_{p-1}^{n-a}$ for $1\leq a\leq n$, and they are
elements of the basis of \fullref{basis}. Denote their dual elements in
$\smash{\smash{\widetilde{\mathcal{K}}_{n,\ast}}^{\ast}}$ by $\wtilde{s}_{n,0}$ and
$c_{n,a}$, respectively (note that $c_{n,n}$ is the unit in $\smash{\widetilde
{\mathcal{K}}_{n,\ast}^{\ast}}$). The topological degree of $e_{2}^{n}$ is
$4(  1+p+\cdots+p^{n-1})  $ and that of $\smash{e_{0}^{a}e_{p-1}^{n-a}}$ is
$2(  p^{n}-p^{a})  .$ It is easy to see, as in the preceding
proposition, that $\smash{\smash{\widetilde{\mathcal{K}}_{n,\ast}}^{\ast}}$ is a polynomial
algebra on the elements $\wtilde{s}_{n,0}$ and $c_{n,a}$, $1\leq a\leq
n-1.$ 

We shall determine the resulting action of the Steenrod algebra. Calculating
with the Nishida formulas in $\smash{\widetilde{\mathcal{K}}_{n,\ast}}$, we obtain,
for $n\geq1$ and $1\leq i\leq n-1$,%
\begin{align*}
d_{p^{n-1}+2p^{n-2}+\cdots+2p+2}\ast e_{2}^{n-1}e_{p+1}  &  =-2e_{2}^{n},\\
d_{p^{n}-p^{i}-p^{i-1}}\ast e_{0}^{i-1}e_{p-1}^{n-i+1}  &  =-e_{0}^{i}%
e_{p-1}^{n-i},\\
d_{p^{n}-p^{n-1}-p^{i}}\ast e_{0}^{i}e_{p-1}^{n-i-1}e_{2p-2}  &  =e_{0}%
^{i}e_{p-1}^{n-i}.
\end{align*}
Converting these to $\smash{\smash{\widetilde{\mathcal{K}}_{n,\ast}}^{\ast}}$ and moving from
the action of $\mathcal{K}$ to that of $\mathcal{A}$, by freely using earlier
lemmas and theorems about $\smash{\widetilde{\mathcal{K}}_{n,\ast}}$, and the
sparseness of $\smash{\smash{\widetilde{\mathcal{K}}_{n,\ast}}^{\ast}}$ in low degrees, we get
the following formulas, for $n\geq1$ and $1\leq i\leq n-1$:%
\begin{align*}
P^{p^{n-1}}\wtilde{s}_{n,0}  &  =2\wtilde{s}_{n,0}c_{n,n-1},\\
P^{p^{i-1}}c_{n,i}  &  =c_{n,i-1}\\
P^{p^{n-1}}c_{n,i}  &  =-c_{n,i}c_{n,n-1}.
\end{align*}
(Note that the occurrence in these formulas of $c_{n,0}$ really denotes
$(  e_{p-1}^{n})  ^{\ast}=\wtilde{s}_{n,0}^{\upnfrac{p-1}{2}}$.)
All other operations $P^{p^{i}}$ are zero on the generators since their images
lie in degrees in which there are no nonzero elements.

These calculated values from our Nishida action on the coalgebra
$\smash{\widetilde{\mathcal{K}}_{n,\ast}}$ coincide with the structure over the
Steenrod algebra of $\mathbb{F}_{p}[t_{1},\ldots,t_{n}]^{\tiny{\smash{\widetilde{SL}_{n}}}}$
as a ring of invariants, as can easily be deduced from the known results in
Kech \cite[page 945]{Kech1993}, \cite[page 280]{Kech1994} and
Wilkerson \cite{Wilkerson}, so they are
isomorphic. Under the isomorphism, $c_{n,a}\in\widetilde{\mathcal{K}}_{n,\ast
}$ maps to the Dickson invariant of the same name in $\smash{\mathbb{F}_{p}%
[t_{1},\ldots,t_{n}]^{\tiny{\smash{\widetilde{SL}_{n}}}}}$, and $\wtilde{s}_{n,0}$ maps to
the remaining polynomial generator of the ring of invariants, which has the
formula $\prod_{i=1}^{n}\widetilde{V}_{i}$ in terms of the elements in the
proof above for the $\widetilde{T}_{n}$ invariants \cite[page
945]{Kech1993}, \cite[page 280]{Kech1994}, see also M\`ui \cite{Mui} and
Wilkerson \cite{Wilkerson}.
\end{proof}

\begin{proof}
[Proof of \fullref{commute}]We will use notation from the two proofs
above. Observe, using our various results, that

\begin{enumerate}
\item In $\smash{\widetilde{\mathcal{K}}_{n,\ast}}$, one has $e_{i}e_{j}=\smash{e_{2}^{2}}$ if
and only if $i=j=2$, and $\smash{e_{i}e_{j}=e_{0}e_{p^{r}(p-1)}}$ if and only if
$(i,j)=(0,p^{r}(p-1))$ or $(i,j)=(p^{r+1}(p-1),0)$ (use the Verschiebung map
$V$, described at the end of the first section of the paper),

\item In $\widetilde{\mathcal{U}}_{n,\ast}^{\ast}$, we have
\[
(  e_{i_{1}}\cdots e_{i_{n}})  ^{\ast}\cdot (  e_{j_{1}}\cdots
e_{j_{n}})  ^{\ast}=(  e_{i_{1}+j_{1}}\cdots e_{i_{n}+j_{n}%
})  ^{\ast},
\]
where all $i$'s and $j$'s are even.
\end{enumerate}

Now, from the first observation we can calculate that
\[
\sigma(c_{n,i})=\sum_{\substack{0\leq j_{1}<\cdots<j_{n-i}\leq n\\j_{0}%
=0}}\left\{  \left(  \prod_{s=1}^{n-i}e_{0}^{j_{s}-j_{s-1}-1}%
e_{(p-1)p^{i+s-j_{s}}}\right)  e_{0}^{n-j_{n-i}}\right\}  ^{\ast}.
\]
And from the second, we find that
\[
\left\{  \left(  \prod_{s=1}^{n-i}e_{0}^{j_{s}-j_{s-1}-1}e_{(p-1)p^{i+s-j_{s}%
}}\right)  e_{0}^{n-j_{n-i}}\right\}  ^{\ast}=\prod_{s=1}^{n-i}\left(
\wtilde{v}_{n,j_{s}-1}\right)  ^{\frac{p-1}{2}p^{^{i+s-j_{s}}}}.
\]
But by the known formulas for the Dickson invariants \cite[page
280]{Kech1994}, \cite[page 224]{Kech}, \cite{Mui},
\[
\tau(c_{n,i})=\sum_{1\leq j_{1}<\cdots<j_{n-i}\leq n}\left\{  \prod
_{s=1}^{n-i}\widetilde{V}_{j_{s}}^{\frac{p-1}{2}p^{^{i+s-j_{s}}}}\right\}  .
\]
A similar calculation shows that $\omega (  \sigma (  \wtilde
{s}_{n,0}))  =\tau(  \wtilde{s}_{n,0})  $, and
the proposition follows.
\end{proof}

\bibliographystyle{gtart}
\bibliography{link}

\end{document}